\documentclass[12pt,a4paper,reqno]{amsart}
\usepackage{amsmath}
\usepackage{amssymb}
\usepackage{amscd}
\usepackage{amsxtra}
\usepackage{hyperref} 
\usepackage{color}
\numberwithin{equation}{section}

\newtheorem{theorem}{Theorem}[section]
\newtheorem{proposition}[theorem]{Proposition}

\newtheorem{lemma}[theorem]{Lemma}

\theoremstyle{definition}

\newtheorem{remark}[theorem]{Remark}

\newcommand{\C}{{\mathbb C}}

\newcommand{\R}{{\mathbb R}}
\newcommand{\Z}{{\mathbb Z}}

\newcommand{\E}{{\mathbb E}}
\newcommand{\Proj}{{\mathbb P}}

\newcommand{\bt}{{\boldsymbol t}}
\newcommand{\mL}{{\mathcal{L}}}
\newcommand{\mO}{{\mathcal{O}}}
\newcommand{\mN}{{\mathcal{N}}}
\newcommand{\mF}{{\mathcal{F}}}

\newcommand{\ev}{\operatorname{ev}}
\newcommand{\Aut}{\operatorname{Aut}}
\newcommand{\ch}{\operatorname{ch}}
\newcommand{\val}{\operatorname{val}}
\newcommand{\virt}{{\rm virt}}
\newcommand{\mov}{{\rm mov}}
\newcommand{\pt}{{\rm pt}}

\newcommand{\Lie}{\operatorname{Lie}}
\newcommand{\Hom}{\operatorname{Hom}}
\newcommand{\Ext}{\operatorname{Ext}}
\newcommand{\Eff}{\operatorname{Eff}}  

\def\ov#1{\overline{#1}}
\def\corr#1{\left\langle{#1}\right\rangle}
\def\parfrac#1#2{{\frac{\partial #1}{\partial #2}}}

\begin{document}

\title{Convergence of Gromov-Witten potentials by localization}
\author{Hiroshi Iritani}
\address{Department of Mathematics, Graduate School of Science, 
Kyoto University, Oiwake-cho, Kitashirakawa, Sakyo-ku, Kyoto, 606-8502, Japan.}
\email{iritani@math.kyoto-u.ac.jp}
\begin{abstract}
We prove that the equivariant genus-$g$ descendant Gromov-Witten potential of a smooth projective variety $X$ 
is a convergent power series in infinitely many variables for each $g$, provided that $X$ is equipped with a torus action with only isolated fixed points and isolated $1$-dimensional orbits. 
\end{abstract}
\maketitle 

\section{Introduction}
Let $X$ be a smooth projective variety equipped with an algebraic action of $T=(\C^\times)^r$. Let $\ov{M}_{g,n}(X,d)$ be the moduli space of genus-$g$, $n$-pointed stable maps to $X$ of degree $d\in H_2(X,\Z)$. This is equipped with the induced $T$-action. 
We have the evaluation map $\ev_i\colon \ov{M}_{g,n}(X,d)\rightarrow X$ at the $i$th marked point 
and the $\psi$-class $\psi_i = c_1^T(\mL_i)$ on $\ov{M}_{g,n}(X,d)$, where $\mL_i \rightarrow \ov{M}_{g,n}(X,d)$ is the line bundle whose fibre at a stable map $f\colon (C;x_1,\dots,x_n) \to X$ is $T^*_{x_i}C$. 
The \emph{equivariant descendant Gromov-Witten invariants} of $X$ are the equivariant integrals 
\[
\corr{\tau_{k_1}\phi_1,\tau_{k_2}\phi_2,\dots,\tau_{k_n}\phi_n}^T_{g,n,d}
=\int_{[\ov{M}_{g,n}(X,d)]_{\virt}}^T 
\prod_{i=1}^n \ev_i^*(\phi_i) \psi_i^{k_i}  
\]
where $\phi_1,\dots,\phi_n\in H^*_T(X,\C)$, $k_1,\dots,k_n\ge 0$ and $[\ov{M}_{g,n}(X,d)]_{\virt}$ is the virtual fundamental class of $\ov{M}_{g,n}(X,d)$. They are elements of $H^*_T(\pt,\C) = \C[\mathfrak{t}_\C]$, i.e.~polynomials on the Lie algebra $\mathfrak{t}_\C = \C^r$. 
The \emph{equivariant genus-$g$ descendant Gromov-Witten potential} $\mF_g$ is a 
generating function of these invariants, given by  
\[
\mF_g(\bt, a;Q) = \sum_{n\ge 0} \sum_{d\in H_2(X,\Z)} 
\corr{\sum_{k=0}^\infty \tau_k \bt_k,\dots,\sum_{k=0}^\infty \tau_k \bt_k}_{g,n,d}^T
\frac{Q^d}{n!},
\]
where $\bt = \{\bt_k\}_{k=0}^\infty$ is an infinite sequence of cohomology classes $\bt_k \in H^*_T(X,\C)$, $a\in \mathfrak{t}_\C$ is the equivariant parameter, and $Q^d$ denotes the element 
of the group ring of $H_2(X,\Z)$ associated with $d\in H_2(X,\Z)$. 
We can think of $\mF_g$ as a formal function on the infinite product $\prod_{k=0}^\infty H^*_T(X,\C)$. 

When $X$ is an algebraic GKM manifold, i.e.~it has only isolated $T$-fixed points and isolated one-dimensional $T$-orbits, we can calculate the Gromov-Witten invariants of $X$ via virtual localization  \cite{kontsevich-enumeration, graber-pandharipande}. 
In this paper, we obtain the following estimate for Gromov-Witten invariants using this technique: 

%\textcolor{red}{Unify the equivariant/nonequivariant statements}
\begin{theorem}\label{thm:main}
Let $X$ be a smooth projective variety admitting an action of an algebraic torus $T$ whose fixed points and $1$-dimensional orbits are isolated. 
% \begin{enumerate}
%\item[(1)] There exists a constant $C(g,X)>0$ depending only on $g$ and $X$ 
%such that for any $\phi_1,\dots,\phi_n\in H^*(X,\C)$,
%\begin{equation}\label{eq:main-estimate}
%\left|\corr{\tau_{k_1}\phi_1,\tau_{k_2}\phi_2,\dots,\tau_{k_n}\phi_n}_{g,n,d}\right|
%\le C(g,X)^{n+|d|+1} n!\, |\phi_1|\cdots |\phi_n|,
%\end{equation}
%where $|d|$ and $|\phi_i|$ are suitable norms on $H_2(X,\R)$ and $H^*(X,\C)$. 
%Consequently, there exists a constant $r=r(g,X)>0$ 
%such that $\mF_g(\bt;Q)$ converges absolutely on the infinite-dimensional polydisc 
%\[
%\Delta_r:=\left\{ (\bt,Q)\;\middle|\; 
%|Q|\le r,\ |\bt_k|\le r \ (k=0,1,2,\dots) \right\}. 
%\]
%\item[(2)] 
Choose a homogeneous $H^*_T(\pt,\C)$-basis $\{\gamma_1,\dots,\gamma_N\}$ of $H^*_T(X,\C)$, norms $|\cdot|$ on $\mathfrak{t}_\C$ and $H_2(X,\R)$ and a K\"ahler class $\omega$ on $X$. 
Then the following hold. 
\begin{itemize} 
\item[(1)] 
There exist constants $C,r>0$ depending only on $g$, $X$ and the choices above such that
\begin{equation}\label{eq:main-estimate}
\left|\corr{\tau_{k_1}\gamma_{h_1},\tau_{k_2}\gamma_{h_2},\dots,\tau_{k_n}\gamma_{h_n}}_{g,n,d}^T(a)\right|
\le C^{n+|d|+1} n! 
\end{equation}
holds whenever $|a|\le r$. 
\item[(2)] Expand $\bt_k = \sum_i t_{k,i} \gamma_i$ and regard $\mF_g(\bt,a;Q)$ as a formal power series in $\{t_{k,i}\}$, the equivariant parameter $a$ and $Q$. 
There exists a constant $r>0$ depending only on $g$, $X$ and the choices above such that 
$\mF_g(\bt,a;Q)$ converges absolutely on the infinite-dimensional polydisc: 
\[
\left\{ (\bt,a,Q)\;\middle|\; 
|a|\le r,\ |Q^d|\le r^{\langle \omega, d\rangle} \ (\text{$d$: effective}),\ |t_{k,i}|\le r \ ( k\ge 0,\ 1\le i\le N) \right\}.
\]
In particular, both $T$-equivariant and non-equivariant quantum cohomology rings of $X$ have convergent structure constants. 
\end{itemize} 
\end{theorem} 

The proof of the theorem is elementary and combinatorial. The virtual localization formula expresses Gromov-Witten invariants as a sum over decorated graphs, and the problem reduces to bounding the number of such graphs. 

\begin{remark} 
(1) Recently, Milanov \cite{Milanov:GKM} proved the convergence of equivariant quantum cohomology of GKM manifolds, by a completely different method (though still based on virtual localization). He established the convergence of the big $J$-function by interpreting the Givental-Brown style recursion \cite{Givental:equivariant,Brown} as a fixed-point problem. 

(2) Examples of algebraic GKM spaces include toric varieties and rational homogeneous spaces $G/P$. The convergence of quantum cohomology for these spaces has been proved in \cite{iritani-colef,CCIT:MS} (for toric varieties and orbifolds) and in \cite{LLYZ} (for $G/P$). 

(3) In this paper, we state and prove our results in the algebraic category. Using the explicit localization formula of Hirschi \cite[Theorem 1.4]{Hirschi}, we can prove the convergence of the Gromov-Witten potentials for Hamiltonian GKM manifolds in the symplectic category by exactly the same argument. 

(4) Zinger \cite{Zinger:genus0, Zinger:energybounds} conjectured the following stronger estimates for \emph{non-equivariant} Gromov-Witten invariants: 
\[
|\corr{\tau_{k_1}\gamma_{h_1},\dots,\tau_{k_n}\gamma_{h_n}}_{g,n,d}| 
\le \frac{n!}{k_1!\cdots k_n!} C_{g,X}^{|d|+n}, 
\]
and proved them for complete intersections in $\Proj^n$ with $g=0$ \cite[Theorem 1]{Zinger:genus0} and for $\Proj^n$ itself in all genera \cite[Theorem 1]{Zinger:energybounds}. 

(5) The \emph{NF-convergence} \cite[Definition 7.5]{coates-iritani:convergence} of the non-equivariant genus-$g$ descendant Gromov-Witten potential $\mF_g^{\rm noneq}$ is equivalent to the following estimates
\[
|\corr{\tau_{k_1}\gamma_{h_1},\dots,\tau_{k_n}\gamma_{h_n}}_{g,n,d}|
\le \frac{n!}{k_1!\cdots k_n!} C_1 C_2^{k_1+\cdots+k_n+n+|d|}   
\] 
for some $C_1,C_2>0$. 
Coates and Iritani \cite[Theorem 1.1]{coates-iritani:convergence} established the NF-convergence of $\mF_g^{\rm noneq}$ on a polydisc independent of $g$ under the assumptions of genus-zero convergence and (generic) semisimplicity. 
Moreover, genus-zero convergence implies the NF-convergence of $\mF_0^{\rm noneq}$ \cite[Theorem 7.9]{coates-iritani:convergence}. Hence our result implies that this estimate holds at $g=0$ for our target space $X$; if the non-equivariant quantum cohomology of $X$ is generically semisimple (which is expected for algebraic GKM manifolds), we know that this estimate holds for all $g\ge 0$ with $C_2$ independent of $g$. 

(6) Recently, McLean \cite{mclean} gave a bound for genus-zero romov-Witten invariants for general smooth projective varieties of the form: 
\[
|\corr{\gamma_{h_1},\dots,\gamma_{h_n}}_{0,n,d}| \le C_1 n! (C_2(\langle \omega, d\rangle +1))!.   
\]
\end{remark}

%\textcolor{green}{On the other hand, the total descendant potential $\mZ(\bt,\epsilon;Q):=\exp(\sum_{g=0}^\infty \epsilon^{g-1} \mF_g(\bt;Q))$ 
%is not necessarily a convergent power series in $\bt, Q,\epsilon$
%since the constant $C(g,X)$ in the above theorem grows rapidly 
%(like $(2g)!$) when $g$ goes to infinity.}  
%
%\textcolor{green}{Recently, the total potentials for local toric Calabi-Yau threefolds 
%were exactly calculated by the ``topological vertex". 
%These calculations suggest that 
%the total potential $\mZ$ should be viewed as 
%an asymptotic expansion in $\epsilon$ of some genuine potential 
%which is holomorphic on some angular domain of $\epsilon$; 
%we examine this phenomenon in Section \ref{sec:higher-genus}. }

%Finally, the author touches on a relation to 
%his previous paper \cite{iritani-colef}.  
%There, using the quantum Lefschetz theorem of Coates and Givental \cite{coates-givental}, 
%he proved that the quantum cohomology of a complete intersection $Y$ in $X$ 
%is partially convergent 
%if the quantum cohomology of the ambient space $X$ is convergent. 
%The motivation of this paper is to give a wide class of examples 
%of ambient spaces $X$ whose quantum cohomology is actually convergent.  

\subsection*{History of this paper and the use of an LLM}
The main results of the paper have been announced on several occasions since 2005 (see e.g.~\cite{iritani-colef, coates-iritani:convergence}). The original motivation was to supplement my result \cite{iritani-colef} on the compatibility of convergence with the quantum Lefschetz theorem \cite{coates-givental}. 

Recently, Todor Milanov informed me of his convergence result \cite{Milanov:GKM} and asked me about the status of the announced results. Prompted by this, I provided a large language model (LLM) with my handwritten notes and an unfinished manuscript, asking it to draft a mathematical paper from them. The LLM quickly generated a text following the flow of my notes, adding several improvements and correcting errors in them (consequently, the present result is weaker than what was announced in \cite[Remark 7.7]{coates-iritani:convergence}). In particular, it completed the proof in the equivariant case, which was not fully addressed in my notes, by using the homogeneity of equivariant correlators to pass from an estimate on a ball $B_0$ away from the origin to one on a ball centred at the origin.

I substantially revised the text generated by the LLM and verified its mathematical correctness.

\subsection*{Acknowledgements} 
I thank Todor Milanov for informing me of his work on the convergence of quantum cohomology for GKM manifolds. I also thank my colleagues from my graduate school days for the helpful discussions we had back then. 

\section{Estimate of Hodge integrals}\label{sec:hodge}
The \emph{Hodge integrals} are intersection numbers on the Deligne-Mumford moduli space $\ov{M}_{g,n}$ of genus-$g$, $n$-pointed stable curves, given by 
\[
\corr{\tau_{k_1},\dots,\tau_{k_n};\prod_{i=1}^s \lambda_{l_i}}_{g,n}
=\int_{\ov{M}_{g,n}} \prod_{i=1}^n \psi_i^{k_i} 
\prod_{i=1}^s \lambda_{l_i} 
\]
where $\lambda_l := c_l(\E_{g,n})$ is the $\lambda$-class, the $l$th Chern class of the Hodge bundle $\E_{g,n}\rightarrow \ov{M}_{g,n}$ whose fibre at a curve $C$ equals $H^0(C,\omega_C)$. 
In this section, we establish the following estimate: 
\begin{proposition}\label{prop:hodge-estimate}
For each $g\ge 0$ and $s\ge 0$, there exists a constant $A(g,s)>0$ such that 
for any $n\ge 0$, $k_1,\dots,k_n\ge 0$, and $0\le l_1,\dots,l_s\le g$,
\begin{equation}\label{eq:hodge-bound}
\left|\corr{\tau_{k_1},\dots,\tau_{k_n};\prod_{i=1}^s \lambda_{l_i}}_{g,n}
\right|
\le \frac{n!}{k_1!\cdots k_n! } A(g,s)^{n+1}.
\end{equation}
\end{proposition}

\begin{remark}\label{rem:genus0-elementary}
When $g=0$ and $s=0$, we have the explicit formula 
\[
\corr{\tau_{k_1},\dots,\tau_{k_n}}_{0,n} = \frac{(n-3)!}{k_1!\cdots k_n!}
\]
for $k_1+\cdots+k_n = n-3$ and $n\ge 3$, so \eqref{eq:hodge-bound} holds 
elementarily with $A(0,0)=1$. 
Similarly, for fixed $g\ge 0$ and $s=0$, one can give an elementary proof of 
a weaker bound $\corr{\tau_{k_1},\dots,\tau_{k_n}}_{g,n} \le C(g)\frac{(2g+n)!}{(2g)!}$ 
by induction on $n$ using the string and dilaton equations:
\begin{align*}
\corr{\tau_0\tau_{k_1}\cdots\tau_{k_l}}_{g,l+1} &= \sum_{i=1}^l \corr{\tau_{k_1}\cdots\tau_{k_i-1}\cdots\tau_{k_l}}_{g,l}, \\
\corr{\tau_1\tau_{k_1}\cdots\tau_{k_l}}_{g,l+1} &= (2g-2+l)\corr{\tau_{k_1}\cdots\tau_{k_l}}_{g,l},
\end{align*}
which reduces any correlator to the finitely many correlators with $k_i\ge 2$ for all $i$.  
%(where $n\le 3g-3$ since $\sum_{i=1}^n (k_i-1) = 3g-3$). 
However, to obtain the sharp factorial denominator $k_1!\cdots k_n!$ needed for 
our localization sum in higher genus, we use Kontsevich's combinatorial model below.
\end{remark}

\begin{proof}[Proof of Proposition \ref{prop:hodge-estimate}]
First we consider the case where $s=0$ (i.e.\ there are no $\lambda$-classes).
We use Kontsevich's identity \cite{kontsevich-airy}, which calculates
descendant intersection numbers by counting ribbon graphs. 
A {\it ribbon graph} is a connected graph 
such that for each vertex $v$, the set of half-edges incident to $v$ 
is equipped with a cyclic order.  
One can naturally associate a compact oriented surface $\Sigma$ 
to a ribbon graph $G$: 
$G$ is embedded in $\Sigma$ so that
the cyclic order at each vertex corresponds to the orientation 
of $\Sigma$ and 
the complement $\Sigma\setminus G$ is a disjoint union of $2$-cells (discs). 
In other words, $\Sigma$ is obtained from $G$ by attaching $2$-cells along the boundary cycles. 
An \emph{$n$-marked ribbon graph} is a ribbon graph $G$ equipped with a marking $\pi_0(\Sigma \setminus G)\cong \{1,\dots,n\}$. 
Define $G^3_{g,n}$ to be the set of isomorphism classes of trivalent 
$n$-marked ribbon graphs such that the associated oriented surface 
is of genus $g$ 
(note that a graph $G\in G^3_{g,n}$ has $2(2g-2+n)$ vertices and $3(2g-2+n)$ edges).   
Let $t_i$ be a variable corresponding to the $i$th marked cell. 
For an edge $e$ of $G$ adjacent to the $i$th and $j$th cells, 
we write $\tilde{t}(e)=t_i+t_j$. 
Kontsevich's identity \cite{kontsevich-airy} states that 
\[ 
\sum_{k_1+\cdots+k_n=3g-3+n} 
\corr{\tau_{k_1},\dots,\tau_{k_n}}_{g,n} \prod_{i=1}^n 
\frac{(2k_i-1)!!}{t_i^{2k_i+1}} 
=\sum_{G\in G_{g,n}^3} \frac{2^{2g-2+n}}{\sharp \Aut(G)} 
\prod_{e\in E(G)} \frac{1}{\tilde{t}(e)}.
\]
Here, $\Aut(G)$ is the automorphism group of $G$ 
as a marked ribbon graph, $E(G)$ is the set of edges of $G$ 
and $(2k-1)!!:=(2k)!/(2^k k!)$ (with $(-1)!!:=1$). 
By Arakelov's theorem \cite{arakelov}, extended by Mumford \cite[Theorem 4.1]{mumford-enumerative} and Arbarello-Cornalba-Griffiths \cite[Chapter XIV, Corollary 5.14]{acgII}, 
$\psi_i$ is nef on $\ov{M}_{g,n}$ and hence $\corr{\tau_{k_1},\dots,\tau_{k_n}}_{g,n}\ge 0$.  
Therefore, by putting $t_i=1$ for all $i=1,\dots,n$ in Kontsevich's formula and using 
$(2k-1)!!\ge k!$, we obtain 
\begin{align*}
0\le \corr{\tau_{k_1},\dots,\tau_{k_n}}_{g,n} 
& \le \frac{1}{\prod_{i=1}^n (2k_i-1)!!} 
\sum_{G\in G^3_{g,n}} \frac{2^{2g-2+n}}{\sharp \Aut(G)}
\frac{1}{2^{3(2g-2+n)}} \\
& \le \frac{\sharp G^3_{g,n}}{k_1!\cdots k_n!} 2^{-2(2g-2+n)}.
\end{align*}
On the other hand, 
Zhi-Cheng Gao \cite{gao} studied the number $T_g(n)$ 
of rooted triangular maps with $n$ vertices 
on a genus $g$ oriented surface. 
He showed the asymptotics 
\[
T_g(n) \sim t_g n^{5(g-1)/2}(12\sqrt{3})^n,  
\]
as $n \to \infty$ for fixed $g$, 
where $t_g>0$ is a constant depending only on $g$. 
Here a triangular map on a surface is 
a connected graph embedded in the surface in 
such a way that each face (connected component of the complement) is a triangle 
($3$-gon), where loops and multiple edges are allowed. 
A triangular map is rooted if it is given a distinguished 
edge, a direction along the edge, and a side of the edge. 
One obtains a trivalent ribbon graph from a triangular map 
by taking its dual graph. 
Because the number of vertices of a triangular map equals  
the number of $2$-cells of the corresponding dual ribbon graph, 
and 
%choosing a root bounds the automorphism factor while 
choosing a labelling of the 
$n$ vertices gives at most $n!$ choices, we have the inequality 
\[
\sharp G_{g,n}^3 \le n! T_g(n).
\]
The conclusion \eqref{eq:hodge-bound} for $s=0$ follows immediately from 
these inequalities. 

Next we consider the case where $s>0$. 
We use the differential equations of Faber and Pandharipande \cite{faber-pandharipande}, 
which determine the correlators containing $\lambda$-classes 
inductively from the $s=0$ case. 
Define a generating function $F_g$ by  
\[
F_g(\{t_k\},\{u_{2m-1}\}) 
=\sum_{n\ge 0} 
\corr{\sum_{k=0}^\infty \tau_k t_k, \dots, \sum_{k=0}^\infty \tau_k t_k
      ;\exp\left(\sum_{m\ge 1} u_{2m-1}\ch_{2m-1}(\E_{g,n})\right)}
_{g,n}\frac{1}{n!}.  
\]
Note that $\ch_{2m}(\E_{g,n})=0$ for $m\ge 1$ by Mumford's relation $c(\E_{g,n}) c(\E_{g,n}^\vee)=1$. 
Using Mumford's Grothendieck-Riemann-Roch calculation \cite{mumford-enumerative}, 
Faber and Pandharipande \cite{faber-pandharipande} showed that $F_g$ satisfies 
the differential equations  
\begin{align}\label{eq:FP-pde}
\parfrac{F_g}{u_{2l-1}} &= 
\frac{B_{2l}}{(2l)!}
\Bigg[ \parfrac{F_g}{t_{2l}}- \sum_{i\ge 0} t_i \parfrac{F_g}{t_{i+2l-1}} \nonumber \\
&\quad +\frac{1}{2} \sum_{i=0}^{2l-2}(-1)^i 
\left(\frac{\partial^2 F_{g-1}}{\partial t_i \partial t_{2l-2-i}}
+\sum_{g=g_1+g_2} \parfrac{F_{g_1}}{t_i}\parfrac{F_{g_2}}{t_{2l-2-i}}
\right)\Bigg]
\end{align}
for all $l\ge 1$, where $B_{2l}$ are the Bernoulli numbers. 
Since the $l$th Chern class $c_l$ 
is a universal polynomial in $\ch_1,\dots,\ch_l$ of weighted degree $l$ 
and $\lambda_l=0$ for $l>g$, 
it suffices to show that there exists a constant $\widetilde{B}(g,s)>0$ such that 
\begin{equation}\label{eq:ch-bound}
\left|\corr{\tau_{k_1},\dots,\tau_{k_n};
\prod_{j=1}^s \ch_{2l_j-1}(\E_{g,n})}_{g,n}\right|
\le \frac{n!}{k_1!\cdots k_n!}\widetilde{B}(g,s)^{n+1}
\end{equation}
for any $1\le l_1,\dots,l_s \le (g+1)/2$. 
The left-hand side of \eqref{eq:ch-bound} equals the absolute value of 
the coefficient of $t_{k_1}\cdots t_{k_n}/n!$ in 
$\partial^s F_{g}/\partial u_{2l_1-1}\cdots \partial u_{2l_s-1}|_{u=0}$. 
Differentiating \eqref{eq:FP-pde} with respect to $u_{2l_2-1},\dots,u_{2l_s-1}$ at $u=0$, we can express $\partial^s F_{g}/\partial u_{2l_1-1}\cdots \partial u_{2l_s-1}|_{u=0}$ in terms of 
derivatives of $F_{g'}$ ($g'\le g$) with at most $s-1$ derivatives in the $u$-variables 
(using the Leibniz rule on the quadratic term). 
Therefore, the estimate \eqref{eq:ch-bound} follows by induction on $s\ge 0$ from 
the differential equation \eqref{eq:FP-pde} and Lemma \ref{lem:diff-op-bound} below. 
\end{proof}

\begin{lemma}\label{lem:diff-op-bound}
Let $F(t_0,t_1,\dots)$ and $G(t_0,t_1,\dots)$ be formal power series of the form 
\[
F(t_0,t_1,\dots) = \sum_{n\ge 0} \frac{1}{n!} \sum_{k_1,\dots,k_n\ge 0} C_{k_1,\dots,k_n}^F t_{k_1}\cdots t_{k_n}
\]
(and similarly for $G$ with coefficients $C_{k_1,\dots,k_n}^G$),
where $C_{k_1,\dots,k_n}^F$ and $C_{k_1,\dots,k_n}^G$ are symmetric tensors satisfying 
\[
|C_{k_1,\dots,k_n}^F|\le A^{n+1} \frac{n!}{k_1!\cdots k_n!}, \qquad 
|C_{k_1,\dots,k_n}^G|\le A^{n+1} \frac{n!}{k_1!\cdots k_n!}
\]
for some constant $A>1$. 
Then for any integers $a,b\ge 0$ and $j\ge 1$, the symmetric coefficients $C'_{k_1,\dots,k_n}$ 
of the power series 
\[
\parfrac{F}{t_a}, \quad \sum_{i=0}^\infty t_i\parfrac{F}{t_{i+j}}, \quad 
\frac{\partial^2 F}{\partial t_a\partial t_b}, \quad 
\parfrac{F}{t_a}\parfrac{G}{t_b}
\]
satisfy the estimate 
\[
|C'_{k_1,\dots,k_n}|\le (A')^{n+1} \frac{n!}{k_1!\cdots k_n!}
\]
for a constant $A'>1$ depending only on $A$.   
\end{lemma}
\begin{proof}
The power series $F$ can be rewritten in multi-index form as 
\[
F(t_0,t_1,\dots)=\sum_{m_0,m_1,\dots\ge 0} C_F(m_0,m_1,\dots )
\frac{t_0^{m_0}}{m_0!} \frac{t_1^{m_1}}{m_1!}\cdots,  
\]
where  
$C_F(m_0,m_1,\dots)
:=C^F_{\underbrace{\scriptstyle 0,\dots,0}_{m_0},\underbrace{\scriptstyle 1,\dots,1}_{m_1},2,\dots}$. 
The hypothesis on $C^F_{k_1,\dots,k_n}$ is equivalent to the estimate 
\begin{equation}\label{eq:multi-index-bound}
|C_F(m_0,m_1,\dots)|\le \frac{(\sum_{i\ge 0} m_i)!}{\prod_{i\ge 0}(i!)^{m_i}} A^{1+\sum_{i\ge 0} m_i}.
\end{equation}

(1) For $\partial F/\partial t_a$, we have 
\[
\parfrac{F}{t_a}=\sum_{m_0,m_1,\dots} C_F(m_0,m_1,\dots) \frac{t_0^{m_0}}{m_0!}
\cdots  \frac{t_a^{m_a-1}}{(m_a-1)!}\cdots. 
\]
The coefficient of $\prod_{i\ge 0} \frac{t_i^{m_i}}{m_i!}$ 
in $\partial F/\partial t_a$ is $C_F(m_0,\dots,m_a+1,\dots)$, which satisfies
\[
|C_F(m_0,\dots,m_a+1,\dots)| \le 
\frac{(1+\sum m_i)!}{a!\prod_{i\ge 0} (i!)^{m_i}} A^{2+\sum m_i}
\le \frac{(\sum m_i)!}{\prod_{i\ge 0} (i!)^{m_i}} {(eA^2)}^{1+\sum m_i},
\] 
where we used $1+n \le e^n$ for $n=\sum m_i \ge 0$ and $a!\ge 1$. 
Thus we can take $A'=eA^2$ for $\partial F/\partial t_a$. 

(2) Applying (1) twice shows that $\partial^2 F/\partial t_a \partial t_b$ satisfies the bound 
with $A' = e(eA^2)^2 = e^3 A^4$.

(3) For $\sum_{i=0}^\infty t_i \partial F/\partial t_{i+j}$ ($j\ge 1$), the coefficient of 
$\prod_{r\ge 0} \frac{t_r^{m_r}}{m_r!}$ is 
\[
\sum_{i\ge 0} m_i C_F(m_0,\dots,m_i-1,\dots,m_{i+j}+1,\dots).
\]
Using $i!/(i+j)! \le 1$ for $j\ge 1$, its absolute value is bounded by 
\begin{align*}
\sum_{i\ge 0} m_i \frac{(\sum_r m_r)!}{\prod_r (r!)^{m_r}} \frac{i!}{(i+j)!} A^{1+\sum_r m_r}
&\le \left(\sum_{i\ge 0} m_i\right) \frac{(\sum_r m_r)!}{\prod_r (r!)^{m_r}} A^{1+\sum_r m_r} \\
&\le \frac{(\sum_r m_r)!}{\prod_r (r!)^{m_r}} (eA)^{1+\sum_r m_r}.
\end{align*}

(4) Finally, by (1), to bound $(\partial F/\partial t_a)(\partial G/\partial t_b)$ it suffices to bound the product $F\cdot G$. 
The coefficient of $\prod_{i\ge 0} \frac{t_i^{m_i}}{m_i!}$ in $F\cdot G$ is 
\[
\sum_{m_i = k_i + l_i} \left(\prod_{i\ge 0} \binom{m_i}{k_i}\right) C_F(k_0,k_1,\dots) C_G(l_0,l_1,\dots).
\]
Using \eqref{eq:multi-index-bound} for $F$ and $G$ together with $\binom{\sum m_i}{\sum k_i} \ge 1$ and $\sum_{k_i\le m_i} \prod_i \binom{m_i}{k_i} = 2^{\sum m_i}$, we obtain 
\begin{align*}
\left|\sum_{m_i = k_i + l_i} \left(\prod_{i\ge 0} \binom{m_i}{k_i}\right) C_F(k) C_G(l)\right|
&\le \sum_{m_i = k_i + l_i} \left(\prod_{i\ge 0} \binom{m_i}{k_i}\right) \frac{(\sum k_i)!(\sum l_i)!}{\prod_i (i!)^{m_i}} A^{2+\sum m_i} \\
&= \frac{(\sum m_i)!}{\prod_i (i!)^{m_i}} A^{2+\sum m_i} \sum_{m_i = k_i + l_i} \frac{\prod_i \binom{m_i}{k_i}}{\binom{\sum m_i}{\sum k_i}} \\
&\le \frac{(\sum m_i)!}{\prod_i (i!)^{m_i}} (2A^2)^{1+\sum m_i}.
\end{align*}
This completes the proof of the lemma.
\end{proof}

\section{Localization and sums over graphs}\label{sec:localization}
In this section, we apply the virtual localization formula of Graber and Pandharipande \cite{graber-pandharipande} 
(following Kontsevich \cite{kontsevich-enumeration} in genus $0$) to prove Theorem \ref{thm:main}.

\subsection{Torus action and fixed loci of the moduli space}
Let $X$ be a smooth projective variety of complex dimension $s = \dim_\C X$ equipped with an algebraic action 
of the torus $T = (\C^\times)^r$ %(or the corresponding compact torus $T_\R = (S^1)^r$) 
such that the fixed point locus $X^T$ and the $1$-dimensional $T$-orbits are isolated. 
Then $X^T$ is a finite set, and the set $E(X)$ of closures of $1$-dimensional $T$-orbits in $X$ is also finite. 
Each $l \in E(X)$ is a smooth $T$-invariant rational curve ($l \cong \Proj^1$) containing exactly two fixed points 
$0_l, \infty_l \in X^T$. 

Let $\mathfrak{t}_\C = \Lie(T)$ and let $\mathfrak{t}^*_\Z = \Hom(T, \C^\times) \subset \mathfrak{t}_\C^*$ denote the weight lattice of $T$. 
At each fixed point $p\in X^T$, the tangent space $T_p X$ decomposes into $1$-dimensional weight spaces 
with weights $w_{p,1},\dots,w_{p,s} \in \mathfrak{t}^*_\Z \setminus \{0\}$. 
Note that any two weights $w_{p,i}$ and $w_{p,j}$ ($i\neq j$) at the same fixed point $p$ 
are linearly independent in $\mathfrak{t}_\C^*$ 
(otherwise, there would be a continuous family of $1$-dimensional $T$-orbits). 

For each $1$-dimensional orbit closure $l\in E(X)$, choose homogeneous coordinates $[z_1:z_2]$ on $l\cong \Proj^1$ 
such that $0_l = [1:0]$ and $\infty_l = [0:1]$. 
Let $w(l) \in \mathfrak{t}^*_\Z \setminus \{0\}$ be the weight of the tangent space $T_{0_l} l$, so that the $T$-action on $l$ is given by 
\[
t\cdot [z_1:z_2] = [z_1 : t^{w(l)} z_2], \qquad t\in T.
\]
The normal bundle $\mN_{l/X}$ of $l$ in $X$ splits equivariantly into a direct sum of $T$-equivariant line bundles:
\[
\mN_{l/X} \cong \bigoplus_{i=1}^{s-1} \mL_l^{(i)}, \qquad \mL_l^{(i)} \cong \mO_{\Proj^1}(m(l,i)) \quad (m(l,i)\in \Z).
\]
Explicitly, $\mL_l^{(i)}$ is the quotient of $(\C^2\setminus\{0\})\times \C$ by $(z_1,z_2,v)\sim (\zeta z_1, \zeta z_2, \zeta^{m(l,i)}v)$ ($\zeta\in\C^\times$) 
equipped with the $T$-action 
\[
t\cdot [z_1, z_2, v] = [z_1, t^{w(l)} z_2, t^{w(l,i)} v],
\]
where $w(l,i)\in \mathfrak{t}^*_\Z$ is the weight of the fibre $(\mL_l^{(i)})_{0_l}$ at $0_l$. 
Then the weight decompositions of $TX$ at $0_l$ and $\infty_l$ are given by
\[
T_{0_l}X = \C_{w(l)} \oplus \bigoplus_{i=1}^{s-1} \C_{w(l,i)}, \qquad 
T_{\infty_l}X = \C_{-w(l)} \oplus \bigoplus_{i=1}^{s-1} \C_{w(l,i)-m(l,i)w(l)}
\]
where $\C_w$ denotes the $1$-dimensional $T$-representation of weight $w$. 
Since the tangent weights at $0_l$ are pairwise linearly independent, $w(l)$ and $w(l,i)$ are linearly independent in $\mathfrak{t}_\C^*$ 
for every $l\in E(X)$ and $1\le i\le s-1$.

Now consider the induced $T$-action on the moduli space $\ov{M}_{g,n}(X,d)$. 
The connected components of the fixed locus $\ov{M}_{g,n}(X,d)^T$ are labelled by isomorphism classes of decorated graphs 
$\hat{\Gamma} = (\Gamma, p, l, d, g, m)$ defined as follows:
\begin{enumerate}
\item $\Gamma$ is a connected graph with vertex set $V(\Gamma)$, edge set $E(\Gamma)$, and tail (half-edge) set $T(\Gamma)$ with $\sharp T(\Gamma)=n$. We allow $\Gamma$ to have multiple edges but no edge loops. A \emph{flag} of $\Gamma$ is a pair $(v,e)$ of a vertex $v\in V(\Gamma)$ and an edge $e\in E(\Gamma)$ incident to $v$; we write $F(\Gamma)$ for the set of flags of $\Gamma$. 
For each vertex $v\in V(\Gamma)$, let $F_v$ denote the set of flags $F=(v,e)\in F(\Gamma)$ incident to $v$, $\val(v) := \sharp F_v$, 
let $T_v \subset T(\Gamma)$ be the set of tails attached to $v$, and $n(v) := \sharp T_v + \val(v)$. 

\item $p\colon V(\Gamma)\to X^T$ assigns a fixed point $p(v)\in X^T$ to each vertex $v$ 
(corresponding to a maximal connected subset $C_v\subset C$ contracted to $p(v)$). 

\item $l\colon E(\Gamma)\to E(X)$ assigns a $1$-dimensional orbit closure $l_e := l(e)\in E(X)$ to each edge $e$ 
(corresponding to a non-contracted irreducible component $C_e\cong \Proj^1$ of $C$) such that if $e$ connects $v_1$ and $v_2$, 
then $l_e$ connects $p(v_1)$ and $p(v_2)$ (in particular, $p(v_1)\neq p(v_2)$). 

\item $d\colon E(\Gamma)\to \Z_{>0}$ assigns a covering degree $d_e := d(e)\ge 1$ to each edge $e$ such that $f|_{C_e}\colon C_e\to l_e$ 
is the unique $T$-equivariant degree $d_e$ cover ramified only over $0_{l_e}$ and $\infty_{l_e}$, and 
\[
\sum_{e\in E(\Gamma)} d_e [l_e] = d \in H_2(X,\Z). 
\]
\item $g\colon V(\Gamma)\to \Z_{\ge 0}$ assigns an arithmetic genus $g(v)$ of $C_v$ to each vertex $v$ such that 
\[
\sum_{v\in V(\Gamma)} g(v) + \dim H^1(\Gamma) = g. 
\]
Here we set $g(v)= 0$ if $C_v$ is a point.   
\item $m\colon \{1,\dots,n\}\xrightarrow{\sim} T(\Gamma)$ is a bijection labelling the tails. We shall identify elements of $T(\Gamma)$ with integers $1\le i\le n$. 
\end{enumerate}

A vertex $v\in V(\Gamma)$ is called {\it stable} if $2g(v)-2+n(v)>0$, in which case $\dim C_v =1$. It is called {\it unstable} if $g(v)=0$ and $n(v)\le 2$, in which case $\dim C_v =0$ and $C_v$ is a single point on $C$ that is either 
\begin{itemize} 
\item a free endpoint of an edge if $\sharp T_v=0, \val(v)=1$, 
\item a marked point on an edge if $\sharp T_v=1,\val(v)=1$, or 
\item a node between two edges if $\sharp T_v=0,\val(v)=2$. 
\end{itemize} 

For each decorated graph $\hat{\Gamma}$, put 
\[
M_{\hat{\Gamma}} := \prod_{\substack{v\in V(\Gamma) \\ 2g(v)-2+n(v)>0}} \ov{M}_{g(v), n(v)}.
\]
There is a canonical finite morphism $\iota_{\hat{\Gamma}}\colon M_{\hat{\Gamma}}\to \ov{M}_{g,n}(X,d)$ whose image is the fixed component 
$F_{\hat{\Gamma}} \cong M_{\hat{\Gamma}} / A_{\hat{\Gamma}}$, where the automorphism group $A_{\hat{\Gamma}}$ fits into the exact sequence 
\[
1 \longrightarrow \prod_{e\in E(\Gamma)} (\Z/d_e\Z) \longrightarrow A_{\hat{\Gamma}} \longrightarrow \Aut(\hat{\Gamma}) \longrightarrow 1.
\]
Here $\Aut(\hat{\Gamma})$ is the automorphism group of $\Gamma$ preserving the decorations $(p,l,d,g,m)$. 
In particular, the order $\sharp A_{\hat{\Gamma}}$ satisfies 
\begin{equation}\label{eq:aut-order}
\sharp A_{\hat{\Gamma}} = \sharp \Aut(\hat{\Gamma}) \prod_{e\in E(\Gamma)} d_e \ge \sharp \Aut(\hat{\Gamma}).
\end{equation}

\subsection{The virtual normal bundle}
We compute the equivariant Euler class $e_T(N_{\hat{\Gamma}}^{\virt})$ of the virtual normal bundle 
$N_{\hat{\Gamma}}^{\virt}$ on $M_{\hat{\Gamma}}$. 
The calculation of $e_T(N_{\hat{\Gamma}}^{\virt})$ has been carried out for $\Proj^s$ by Graber-Pandharipande \cite{graber-pandharipande} and for our target spaces $X$ (algebraic GKM manifolds) by Liu-Sheshmani \cite{liu-sheshmani}. We briefly review the computation.

The tangent space $T^1$ and the obstruction space $T^2$ at a stable map $f\colon (C,x_1,\dots,x_n) \to X$ fit into the exact sequence 
\begin{align*} 
0\to &\Ext^0(\Omega_C(x_1+\cdots+x_n),\mO_C) \to H^0(C,f^*TX) \to T^1 \\ 
\to &\Ext^1(\Omega_C(x_1+\cdots+x_n),\mO_C) \to H^1(C,f^*TX) \to T^2 \to 0.  
\end{align*} 
Let $f$ be a stable map representing a point in the $T$-fixed component $F_{\hat{\Gamma}} = M_{\hat{\Gamma}}/A_{\hat{\Gamma}}$ associated with $\hat{\Gamma}$. The virtual normal bundle $N^{\virt}_{\hat{\Gamma}}$ at $f$ is given by the moving part of $T^1- T^2$ with respect to the $T$-action, i.e.~
%In the representation ring of $T$, the virtual tangent space at $(f\colon (C, x_1,\dots,x_n)\to X)\in M_{\hat{\Gamma}}/A_{\hat{\Gamma}}$ is given by \textcolor{red}{Wrong} 
%\begin{align*}
%T_f \ov{M}_{g,n}(X,d)
%&= [H^0(C, f^*TX)] - [H^1(C, f^*TX)] + \sum_{i=1}^n [T_{x_i}C] \\
%&\quad + \sum_{\substack{y\in C^\alpha \cap C^\beta \\ \alpha \neq \beta}} \left([T_y C^\alpha] + [T_y C^\beta] + [T_y C^\alpha \otimes T_y C^\beta]\right)
%- \sum_{\alpha} [H^0(C^\alpha, TC^\alpha)],
%\end{align*}
%where $C = \bigcup_\alpha C^\alpha$ is the decomposition into irreducible components. 
%The fixed part of $T_f \ov{M}_{g,n}(X,d)$ equals $T M_{\hat{\Gamma}}$, which comes from the deformations of the contracted subsets $C_v$ and special points on them. 
%Thus the virtual normal bundle $N_{\hat{\Gamma},f}^{\virt} = (T_f \ov{M}_{g,n}(X,d))^{\mov}$ at $f$ is given by
\begin{align}
\label{eq:Nvirt-K}
\begin{split} 
N_{\hat{\Gamma},f}^{\virt}
& =  [H^0(C, f^*TX)]^{\mov} - [H^1(C, f^*TX)]^{\mov} \\ 
& \quad - [\Ext^0(\Omega_C(x_1+\cdots+x_n),\mO_C)]^{\mov} + [\Ext^1(\Omega_C(x_1+\cdots+x_n),\mO_C)]^{\mov}. 
\end{split} 
%+ \sum_{x_i \in C_e,\, e\in E(\Gamma)} [T_{x_i} C_e] \nonumber \\
%&\quad + \sum_{F=(v,e)\in F(\Gamma), \, \textcolor{red}{\dim C_v=1}} [T_{x_F} C_e] + \sum_{F=(v,e)\in F(\Gamma),\, \dim C_v=1} [T_{x_F} C_v \otimes T_{x_F} C_e] \\
%&\quad + \sum_{\substack{v\in V(\Gamma),\,\dim C_v=0 \\ E_v = \{e_1, e_2\}, \, F_i:=(v,e_i)}} [T_{x_{F_1}} C_{e_1} \otimes T_{x_{F_2}} C_{e_2}]
%- \sum_{e\in E(\Gamma)} [H^0(C_e, TC_e)]^{\mov}, \nonumber
\end{align}
%where $x_F \in C$ is the unique intersection point of $C_e$ and $C_v$ for a flag $F=(v,e)$. 

To compute $[H^0(C,f^*TX)] - [H^1(C,f^*TX)]$, consider the normalization exact sequence of $C$ along the nodes $x_F = C_v \cap C_e$ 
indexed by flags $F=(v,e)$:
\[
0 \longrightarrow \mO_C \longrightarrow \bigoplus_{v\in V(\Gamma)} \mO_{C_v} \oplus \bigoplus_{e\in E(\Gamma)} \mO_{C_e}
\longrightarrow \bigoplus_{F=(v,e)\in F(\Gamma)} \mO_{x_F} \longrightarrow 0.
\]
Tensoring with $f^*TX$ and taking the long exact sequence in cohomology yields
\begin{multline*} %\label{eq:H0-H1-decomp}
[H^0(C, f^*TX)]  - [H^1(C, f^*TX)] 
 = \sum_{e\in E(\Gamma)} \left([H^0(C_e, f^*TX)] - [H^1(C_e, f^*TX)]\right) \nonumber \\
 - \sum_{v\in V(\Gamma)} (\val(v)-1)[T_{p(v)}X]
- \sum_{\substack{v\in V(\Gamma) \\ \dim C_v = 1}} [\E_{g(v), n(v)}^\vee \otimes T_{p(v)}X],
\end{multline*}
where we used $H^0(C_v, f^*TX) = T_{p(v)}X$ and $H^1(C_v, f^*TX) = H^1(C_v, \mO_{C_v})\otimes T_{p(v)}X \cong \E_{g(v),n(v)}^\vee \otimes T_{p(v)}X$.

For each edge $e\in E(\Gamma)$, write $w_e := w(l_e)$, $w_e^i := w(l_e,i)$, and $m_e^i := m(l_e,i)$. 
In homogeneous coordinates $[z_1:z_2]$ on $C_e\cong \Proj^1$ with $0_e = [1:0]$ and $\infty_e = [0:1]$, 
the map $f|_{C_e}\colon C_e\to l_e$ is $[z_1:z_2]\mapsto [z_1^{d_e}:z_2^{d_e}]$ and the lifted $T$-action on $C_e$ is 
$t\cdot [z_1:z_2] = [z_1 : t^{w_e/d_e} z_2]$. 
Decomposing $f^*TX|_{C_e} \cong f^*Tl_e \oplus \bigoplus_{i=1}^{s-1} f^*\mL_{l_e}^{(i)}$, we have: 
\begin{itemize}
\item Since $f^*Tl_e \cong \mO_{\Proj^1}(2d_e)$, we have $H^1(C_e, f^*Tl_e)=0$ and $H^0(C_e, f^*Tl_e)$ has a basis of sections 
with $T$-weights 
\[
w_e\left(\frac{j}{d_e} - 1\right), \qquad 0\le j\le 2d_e 
\]
where the weight with $j=d_e$ does not contribute to the moving part. 
\item For each $i\in \{1,\dots,s-1\}$ with $m_e^i \ge 0$, we have $f^*\mL_{l_e}^{(i)}\cong \mO_{\Proj^1}(d_e m_e^i)$, so $H^1(C_e, f^*\mL_{l_e}^{(i)})=0$ 
and $H^0(C_e, f^*\mL_{l_e}^{(i)})$ is spanned by monomials $z_1^j z_2^{d_e m_e^i - j}$ ($0\le j\le d_e m_e^i$) with weights in 
\[
I_e^i := \left\{ w_e^i - w_e\left(m_e^i - \frac{j}{d_e}\right) \;\middle|\; 0\le j\le d_e m_e^i \right\}.
\]
\item For each $i\in \{1,\dots,s-1\}$ with $m_e^i < 0$, both $H^0$ and $H^1$ vanish if $d_e m_e^i = -1$, while if $d_e m_e^i \le -2$, 
$H^0(C_e, f^*\mL_{l_e}^{(i)})=0$ and by Serre duality $H^1(C_e, f^*\mL_{l_e}^{(i)}) \cong H^0(C_e, (f^*\mL_{l_e}^{(i)})^\vee \otimes \Omega_{C_e})^\vee$ 
has weights in 
\[
J_e^i := \left\{ w_e^i - w_e\left(m_e^i + \frac{1+j}{d_e}\right) \;\middle|\; 0\le j\le -2 - d_e m_e^i \right\}.
\]
\end{itemize}
Note that neither $I_e^i$ nor $J_e^i$ contains $0$ since $w_e$ and $w_e^i$ are linearly independent. 

Next we consider the curve automorphism and deformation terms, i.e.~the second line in \eqref{eq:Nvirt-K}. 
We have 
\begin{align*} 
[\Ext^0(\Omega_C(x_1+\cdots+x_n),\mO_C)]^\mov 
& = [H^0(C,TC(-x_1-\cdots-x_n))]^\mov 
\\ 
& = \sum_{F=(v,e) \in F(\Gamma), \, \dim C_v=0,\, n(v) = 1} [T_{x_F}C_e] 
\intertext{and} 
[\Ext^1(\Omega_C(x_1+\cdots+x_n),\mO_C)]^\mov 
& =  \sum_{F=(v,e)\in F(\Gamma),\, \dim C_v=1} [T_{x_F} C_v \otimes T_{x_F} C_e] \\
&\quad + \sum_{\substack{v\in V(\Gamma),\,\dim C_v=0 \\ F_v = \{F_1, F_2\}, \, F_i=(v,e_i)}} [T_{x_{F_1}} C_{e_1} \otimes T_{x_{F_2}} C_{e_2}]
\end{align*} 
%\textcolor{green}{The space $H^0(C_e, TC_e)$ is spanned by $\partial_{z}, z\partial_{z}, z^2\partial_{z}$ with weights 
%$\frac{w_e}{d_e}, 0, -\frac{w_e}{d_e}$, where $z=z_2/z_1$ is the inhomogeneous coordinate on $C_e$. 
%The weight-$0$ term cancels the $j=d_e$ (weight-$0$) term in $H^0(C_e, f^*Tl_e)$. 
%The remaining two terms of weights $\pm \frac{w_e}{d_e}$ in $H^0(C_e, TC_e)$ cancel against $[T_{0_e}C_e] = \lambda^{w_e/d_e}$ 
%and $[T_{\infty_e}C_e] = \lambda^{-w_e/d_e}$ whenever $0_e$ or $\infty_e$ is a special point (a node or a marked point); 
%if $0_e$ (resp.\ $\infty_e$) is non-special (i.e.\ an unstable vertex $v$ with $g(v)=0, \sharp T_v=0, \val(v)=1$), 
%the corresponding factor $\frac{w_e}{d_e}$ (resp.\ $-\frac{w_e}{d_e}$) remains in the denominator of $e_T(N_{\hat{\Gamma}}^{\virt})$.
%}

Combining the above computations and setting $w_F := w_e$ if $p(v)=0_{l_e}$ and $w_F := -w_e$ if $p(v)=\infty_{l_e}$ for a flag $F=(v,e)$, we arrive at the explicit formula:
\begin{align}\label{eq:euler-Nvirt}
\begin{split} 
& \frac{1}{e_T(N_{\hat{\Gamma}}^{\virt})}
= \prod_{e\in E(\Gamma)} \left[ (-1)^{d_e} \left(\frac{d_e^{d_e}}{w_e^{d_e} d_e!}\right)^2 
\frac{\prod_{i:\, d_e m_e^i \le -2} \prod_{\alpha\in J_e^i} \alpha}{\prod_{i:\, m_e^i\ge 0} \prod_{\alpha\in I_e^i} \alpha} \right] \\ 
& \quad \times 
\prod_{v\in V(\Gamma)} e_T(T_{p(v)}X)^{\val(v)-1} 
\prod_{\substack{F=(v,e) \in F(\Gamma) \\ \dim C_v=0,\, n(v)=1}} \frac{w_F}{d_e}  
\prod_{\substack{v\in V(\Gamma): \dim C_v=0 \\ F_v = \{F_1, F_2\}, \, F_i = (v,e_i)}} \frac{1}{\frac{w_{F_1}}{d_{e_1}} + \frac{w_{F_2}}{d_{e_2}}}
\\ 
& \quad \times 
\prod_{\substack{v\in V(\Gamma) \\ \dim C_v = 1}} 
\prod_{i=1}^s \left(\sum_{j=0}^{g(v)} (-1)^j (w_{p(v),i})^{g(v)-j} \lambda_j\right)
\prod_{\substack{F=(v,e)\in F(\Gamma) \\ \dim C_v = 1}} \frac{1}{\frac{w_F}{d_e} - \psi_F}
\end{split} 
\end{align}
where $\psi_F$ denotes the $\psi$-class on $M_{\hat{\Gamma}}$ at the marking associated with a flag $F$. 

%\textcolor{red}{Why $e_T(T_{p(v)}X)^{\val(v)-1}$-terms separate? $e_T(T_{p(v)}X) = \prod_i w_{p(v),i}$}

\subsection{Estimates for a fixed graph}
By equivariant formality \cite{gkm}, $H^*_T(X,\C)$ is a free $H^*_T(\pt,\C)$-module of rank $N=\dim_\C H^*(X,\C)$. 
We choose a homogeneous $H^*_T(\pt,\C)$-basis $\{\gamma_1,\dots,\gamma_N\}$ of $H^*_T(X,\C)$. 
%\textcolor{red}{Do we need a splitting?}
The virtual localization formula \cite{graber-pandharipande} gives
\begin{equation}\label{eq:loc-sum}
\corr{\tau_{k_1}\gamma_{h_1},\dots,\tau_{k_n}\gamma_{h_n}}^T_{g,n,d}
= \sum_{\text{$\hat{\Gamma}$: genus $g$, degree $d$, $n$ tails}} \frac{1}{\sharp A_{\hat{\Gamma}}} \int_{M_{\hat{\Gamma}}} \frac{\prod_{i=1}^n \gamma_{h_i}|_{p(v_i)} \psi_i^{k_i}}{e_T(N_{\hat{\Gamma}}^{\virt})},
\end{equation}
where $v_i \in V(\Gamma)$ denotes the vertex to which the $i$th tail is attached (i.e.~$i\in T_{v_i}$). 
Note that if $v_i$ is a stable vertex, $\psi_i$ restricts to the $\psi$-class on $M_{\hat{\Gamma}}$, 
whereas if $v_i$ is an unstable vertex with incident flag $F_i = (v_i,e_i)$, $\psi_i|_{M_{\hat{\Gamma}}} = -w_{F_i}/d_{e_i}$.   
% Since the left-hand side of \eqref{eq:loc-sum} is independent of the equivariant parameter $\lambda \in \mathfrak{t}_\C$ 
% (or if we wish to bound the equivariant correlators at a fixed $\lambda$), we may evaluate the right-hand side at a suitably chosen $\lambda \in \mathfrak{t}_\C$ with $|\lambda|=1$.

In the following discussion, we fix norms $|\cdot|$ on the vector spaces $\mathfrak{t}_\C$ and $H_2(X,\R)$. We shall evaluate the right-hand side of \eqref{eq:loc-sum} on a closed ball $B_0\subset \mathfrak{t}_\C$ away from the origin. 

\begin{lemma}\label{lem:choice-of-ball}
There exist a closed ball $B_0 \subset \mathfrak{t}_\C$ of radius $0<\delta<1$ not containing the origin 
and constants $C_1, C_2, C_3 > 1$, $0<\epsilon<1$ depending only on $X$ such that for all $a\in B_0$, all decorated graphs $\hat{\Gamma}$, 
all $e\in E(\Gamma)$, all $p\in X^T$, all $1\le i\le s-1$ and all $1\le j,k\le s$, we have: 
\begin{align} 
\label{eq:IJ-bounds}
\begin{split}
& \prod_{\alpha \in I_e^i} |\alpha(a)|^{-1} \le C_1 C_2^{d_e} \quad \text{if $m_e^i\ge 0$}, \qquad
\prod_{\alpha \in J_e^i} | \alpha(a) | \le C_1 C_2^{d_e} \quad\text{if $d_e m_e^i \le -2$}. \\ 
& 1/C_3\le |w_e(a)|\le C_3, \quad 1/C_3\le |w_{p,j}(a)|\le C_3, \quad \cos\left(\arg\left(-\tfrac{w_{p,j}(a)}{w_{p,k}(a)}\right)\right) \le 1-\epsilon. 
\end{split} 
\end{align}
\end{lemma}
\begin{proof}
For each $l\in E(X)$ and $1\le i\le s-1$, the weights $w(l)$ and $w(l,i)$ are linearly independent in $\mathfrak{t}^*_\C$. This implies that the quotient $w(l)/w(l,i)$, viewed as a rational function on $\mathfrak{t}_\C$, is not constant. Similarly, for each $1\le j,k\le s$, $-w_{p,k}/w_{p,j}$ is either the constant function $-1$ (which happens when $j=k$) or a nonconstant rational function. 
Therefore we may find a closed ball $B_0\subset \mathfrak{t}_\C$ of radius $\delta>0$ such that 
\begin{itemize} 
\item the linear functions $w(l)$, $w(l,i)$ and $w_{p,j}$ have no zeros on $B_0$; 
\item the rational function $w(l)/w(l,i)$ does not take real values on $B_0$; and 
\item the rational function $-w_{p,k}/w_{p,j}$ does not take values in $\R_{\ge 0}$ on $B_0$ 
\end{itemize} 
for all $l$, $1\le i\le s-1$ and $1\le j,k\le s$. The first and third conditions ensure that there exist $C_3>1$ and $0<\epsilon<1$ satisfying the second line of \eqref{eq:IJ-bounds} (note that $w_e = w(l_e)$).  

Consider the real line segment $I(a) = [w_e^i(a) - m_e^i w_e(a), w_e^i(a)] \subset \C$ for $a\in B_0$; this depends only on $l_e$, $i$ and $a$ (recall that $w_e^i = w(l_e,i)$, $w_e = w(l_e)$ and $m_e^i = m(l_e,i)$).  By the choice of $B_0$, this segment does not pass through $0$. Since $\bigcup_{a\in B_0} I(a)$ is compact and there are only finitely many $l\in E(X)$ and $1\le i\le s-1$, we may find constants $c, C>0$ depending only on $X$ such that 
\[
c\le |z| \le C \quad \text{for all $z\in I(a)$ and $a\in B_0$}.  
\]
This gives uniform bounds $\prod_{\alpha\in I_e^i} |\alpha(a)|^{-1} \le c^{-(d_e m_e^i + 1)}$ (when $m_e^i\ge 0$) and 
$\prod_{\alpha\in J_e^i} |\alpha(a)| \le C^{d_e |m_e^i|-1}$ (when $d_e m_e^i \le -2$) on $B_0$, for all  $\hat{\Gamma}$, $e$ and $i$. They establish the first line of \eqref{eq:IJ-bounds} for some $C_1,C_2>1$. 
\end{proof}

In the rest of Section \ref{sec:localization}, we choose a ball $B_0\subset \mathfrak{t}_\C$ as in Lemma \ref{lem:choice-of-ball}. 
Let $\omega \in H^2(X,\R)$ be a K\"ahler class on $X$. Since every $l\in E(X)$ is an effective curve, 
$c_0 := \min_{l\in E(X)} \langle \omega, [l]\rangle > 0$. 
Choose a constant $C_0>0$ such that $|\langle \omega, d\rangle| \le C_0 |d|$ for all $d\in H_2(X,\R)$. 
Then for any decorated graph $\hat{\Gamma}$ of degree $d = \sum_{e\in E(\Gamma)} d_e [l_e]$, setting $D := \sum_{e\in E(\Gamma)} d_e$, we have 
\begin{equation} 
\label{eq:edge-degree-bound}
\sharp E(\Gamma) \le D = \sum_{e\in E(\Gamma)} d_e \le \frac{1}{c_0} \sum_{e\in E(\Gamma)} d_e \langle \omega, [l_e]\rangle = \frac{\langle \omega, d\rangle}{c_0} \le  C |d|,
\end{equation}
with $C=C_0/c_0$ and hence 
\begin{equation} 
\label{eq:vertex-degree-bound} 
\sharp V(\Gamma) = \sharp E(\Gamma) + 1 - \dim H^1(\Gamma) \le C|d| + 1. 
\end{equation} 

Assuming $a\in B_0$, we estimate each factor of $1/e_T(N_{\hat{\Gamma}}^\virt)$ in \eqref{eq:euler-Nvirt} that does not involve the $\psi$- and $\lambda$-classes (appearing in the first two lines of the right-hand side of \eqref{eq:euler-Nvirt}). The constants $C_3,C_4,\dots$ below depend only on $X$. 
\begin{itemize}
\item Since $\log(d_e! / d_e^{d_e}) \ge \int_1^{d_e} \log x\,dx - d_e \log d_e = 1 - d_e$, 
we have 
\[
\left(\frac{d_e^{d_e}}{| w_e(a)|^{d_e} d_e!}\right)^2 \le e^{-2+2d_e} C_3^{2d_e} 
\le (eC_3)^{2d_e}.
\]
Hence, by \eqref{eq:IJ-bounds}, \eqref{eq:edge-degree-bound}, \eqref{eq:vertex-degree-bound}, 
\begin{align*} 
\left |\prod_{e\in E(\Gamma)} \left[  \left(\frac{d_e^{d_e}}{w_e(a)^{d_e} d_e!}\right)^2 
\frac{\prod_{i:\, d_e m_e^i \le -2} \prod_{\alpha\in J_e^i} \alpha(a)}{\prod_{i:\, m_e^i\ge 0} \prod_{\alpha\in I_e^i} \alpha(a)} \right] \right| 
\le (eC_3)^{2D} C_1^{(s-1)\sharp E(\Gamma)}C_2^{(s-1)D} \le C_4^{|d|+1} 
\end{align*} 
for some $C_4>1$. 
\item For each unstable $2$-valent vertex $v$ ($\dim C_v=0, \sharp T_v=0, \val(v)=2$) with incident flags $F_1=(v,e_1)$ and $F_2=(v,e_2)$, we rewrite 
\[
\frac{1}{\frac{w_{F_1}}{d_{e_1}} + \frac{w_{F_2}}{d_{e_2}}}
= \frac{d_{e_1} + d_{e_2}}{\left(1 + \frac{d_{e_2}}{d_{e_1}}\right) w_{F_1}+ \left(1 + \frac{d_{e_1}}{d_{e_2}}\right)w_{F_2}} 
= \frac{d_{e_1}+d_{e_2}}{t_1 w_{F_1} + t_2 w_{F_2}} 
\]
with $t_1=1+d_{e_2}/d_{e_1}, t_2=1+d_{e_1}/d_{e_2}\ge 1$. By the cosine formula, setting $\theta = \arg(-w_{F_1}(a)/w_{F_2}(a))$, we have 
\begin{align*} 
\left|t_1 w_{F_1}(a) +t_2 w_{F_2}(a)\right|^2 
& = t_1^2 |w_{F_1}(a)|^2 + t_2^2 |w_{F_2}(a)|^2 
-2 t_1t_2|w_{F_1}(a)||w_{F_2}(a)| \cos \theta \\ 
& \ge (t_1 |w_{F_1}(a)|-t_2|w_{F_2}(a)|)^2 + 2\epsilon t_1t_2 |w_{F_1}(a)||w_{F_2}(a)| \\ 
& \ge 2 \epsilon C_3^{-2} \ge \epsilon C_3^{-2}
\end{align*} 
where we used the estimates \eqref{eq:IJ-bounds}.  
Using $d_{e_1}+d_{e_2} \le e^{d_{e_1}+d_{e_2}}$ and $\sum_{v} (d_{e_1}+d_{e_2}) \le 2\sum_{e\in E(\Gamma)} d_e = 2D$, we obtain using \eqref{eq:vertex-degree-bound} 
\[
\prod_{\substack{v:\, \dim C_v=0,\, \sharp T_v=0 \\ F_v=\{F_1,F_2\},\, F_i=(v,e_i)}} \left|\frac{1}{\frac{w_{F_1}(a)}{d_{e_1}} + \frac{w_{F_2}(a)}{d_{e_2}}}\right|
\le e^{2D} \left(\frac{C_3}{\sqrt{\epsilon}}\right)^{\sharp V(\Gamma)}  \le C_5^{|d|+1}
\]
for some $C_5>1$. 
\item Since $\sum_{v\in V(\Gamma)} \val(v) = 2\sharp E(\Gamma)$ and $e_T(T_{p(v)}X) = \prod_{j=1}^s w_{p(v),j}(a)$, we have 
\[
\prod_{v\in V(\Gamma)} |e_T(T_{p(v)}X)|^{\val(v)-1} \prod_{\substack{F=(v,e) \in F(\Gamma) \\ \dim C_v=0,\, n(v)=1}} \left|\frac{w_F(a)}{d_e}  \right| 
 \le C_3^{s(1+2\sharp E(\Gamma))} C_3^{\sharp V(\Gamma)} \le C_6^{|d|+1} 
\]
for some $C_6>1$. 
\end{itemize}
Altogether, whenever $a\in B_0$, we have 
\begin{equation} 
\label{eq:scalar-factor}
\left|\left( 
\renewcommand{\arraystretch}{0.7} 
\begin{array}{c} \tiny \text{the product of the factors} \\
\tiny\text{in the first two lines} \\ 
\tiny \text{of the RHS of \eqref{eq:euler-Nvirt}}\end{array} \right) \right| \le C_7^{1+|d|}
\end{equation} 
for some $C_7>1$.

Finally we estimate the integral over $M_{\hat{\Gamma}}$. We have 
\begin{multline} 
\label{eq:integral-over-MGamma}
\int_{M_{\hat{\Gamma}}} \frac{\prod_{i=1}^n \gamma_{h_i}|_{p(v_i)} \psi_i^{k_i}}{e_T(N_{\hat{\Gamma}}^{\virt})}
= \left( \renewcommand{\arraystretch}{0.7} 
\begin{array}{c} \tiny \text{the product of the factors} \\
\tiny\text{in the first two lines} \\ 
\tiny \text{of the RHS of \eqref{eq:euler-Nvirt}}\end{array} \right) 
\prod_{i=1}^n \gamma_{h_i}|_{p(v_i)} 
\prod_{\substack{i: \dim C_{v_i}=0 \\ F= (v_i,e_i) \in F(\Gamma)}}
\left( - \frac{w_{F}(a)}{d_{e_i}}\right)^{k_i}
\\
\times \prod_{v\in V(\Gamma): \dim C_v=1} 
\left ( \int_{\ov{M}_{g(v),n(v)}} 
\frac{\prod_{i=1}^s \left(\sum_{j=0}^{g(v)} (-1)^j (w_{p(v),i}(a))^{g(v)-j} \lambda_j\right)}
{\prod_{F=(v,e)\in F_v} \left(\frac{w_F(a)}{d_e} - \psi_F\right)}
\prod_{i\in T_v} \psi_i^{k_i} \right).
\end{multline} 
By \eqref{eq:IJ-bounds} and \eqref{eq:scalar-factor}, whenever $a\in B_0$, 
the norm of the first line of \eqref{eq:integral-over-MGamma} can be uniformly bounded by $C_8^{1+n+|d|+\sum_{i=1}^n k_i}$ for some $C_8>1$. 
Let $V_{\rm st}=\{v\in V(\Gamma): \dim C_v=1\}$ be the set of stable vertices and let $F_{\rm st} = \bigcup_{v\in V_{\rm st}} F_v$ be the set of stable flags.    
We expand the second line of \eqref{eq:integral-over-MGamma} as follows: 
\begin{align*} 
\sum_{\substack{j\colon V_{\rm st}\times \{1,\dots,s\}\to \Z_{\ge 0} \\ \text{$j(v,i) \le g(v)$ for all $(v,i)$}}} 
\sum_{b\colon F_{\rm st} \to \Z_{\ge 0}}
\prod_{v\in V_{\rm st}} 
\Biggl( 
\prod_{F=(v,e)\in F_v}\left(\frac{d_{e}}{w_F(a)}\right)^{b(F)+1}\prod_{i=1}^s (-1)^{j(v,i)} (w_{p(v),i}(a))^{g(v)-j(v,i)}  \\ 
\times 
\int_{\ov{M}_{g(v),n(v)}}  \prod_{F\in F_v} \psi_F^{b(F)} 
\prod_{i\in T_v} \psi_i^{k_i} 
\prod_{i=1}^s \lambda_{j(v,i)} 
\Biggr) 
\end{align*}
By Proposition \ref{prop:hodge-estimate}, we have 
\[
\left|\int_{\ov{M}_{g(v), n(v)}} \prod_{F\in F_v} \psi_F^{b(F)} \prod_{i\in T_v}\psi_i^{k_i} \prod_{i=1}^s \lambda_{j(v,i)}\right|
\le A_1(g,s)^{n(v)+1} \frac{n(v)!}{\left(\prod_{i\in T_v} k_i!\right)\left(\prod_{F\in F_v} b(F)!\right)}.
\]
for $A_1(g,s) = \max_{0\le g'\le g}A(g',s)>1$. 
Thus the norm of the second line of \eqref{eq:integral-over-MGamma} is bounded by 
\begin{align*} 
\sum_{\substack{j\colon V_{\rm st}\times \{1,\dots,s\}\to \Z_{\ge 0} \\ \text{$j(v,i) \le g(v)$ for all $(v,i)$}}} 
&\prod_{v\in V_{\rm st}} 
\left( 
C_3^{sg(v)} A_1(g,s)^{n(v)+1} n(v)! 
\prod_{F=(v,e)\in F_v} 
\sum_{b(F)=0}^\infty \frac{1}{b(F)!} \left(\frac{d_e}{|w_F(a)|}\right)^{b(F)+1}
\right) 
\\
&\le (g+1)^{s \sharp V(\Gamma)} C_3^{sg} \prod_{v\in V_{\rm st}} \left( A_1(g,s)^{n(v)+1} n(v)! 
\prod_{F=(v,e) \in F_v} C_3 d_e e^{C_3d_e}\right) \\ 
&\le (g+1)^{s (C|d|+1)} C_3^{sg} A_1(g,s)^{3C|d|+n+1} e^{4C_3 C |d|} \prod_{v\in V(\Gamma)} n(v)!
\end{align*} 
whenever $a\in B_0$, where we used the fact that the number of maps $j\colon V_{\rm st}\times \{1,\dots,s\} \to \Z_{\ge 0}$ is bounded by $(g+1)^{s\sharp V}$, $d_e/|w_F(a)| \le C_3 d_e$ in the second line; \eqref{eq:edge-degree-bound}, \eqref{eq:vertex-degree-bound}, 
$\sum_{v\in V_{\rm st}} (n(v)+1) \le 2 \sharp E(\Gamma) + n + \sharp V(\Gamma) 
\le 3C |d|+n+1$, $C_3 d_e \le e^{C_3d_e}$ and $\sum_{F=(v,e)\in F_{\rm st}} d_e \le 2 D\le 2C|d|$. 

Combining the above estimates, we find that there exists a constant $A_2(g,X)\ge 1$, depending only on $X$ and the genus $g$, such that whenever $a\in B_0$, 
\begin{equation}\label{eq:single-graph-bound}
\left|\int_{M_{\hat{\Gamma}}} \frac{\prod_{i=1}^n \gamma_{h_i}|_{p(v_i)} \psi_i^{k_i}}{e_T(N_{\hat{\Gamma}}^{\virt})}\right|
\le A_2(g,X)^{n+|d|+1 + \sum_{i=1}^n k_i} \prod_{v\in V(\Gamma)} n(v)!.
\end{equation}

%\begin{remark}\label{rem:genus0-multinomial}
%In the genus $0$ case ($g=0$), one can also see \eqref{eq:single-graph-bound} directly from Remark \ref{rem:genus0-elementary}: 
%the integral at each vertex $v$ with $n(v)\ge 3$ equals 
%\[
%\left(\prod_{F\in F_v} \frac{d_e}{w_F}\right) \left(\sum_{F\in F_v} \frac{d_e}{w_F}\right)^{n(v)-3 - \sum_{i\in T_v} k_i} 
%\frac{(n(v)-3)!}{\left(\sum_{F\in F_v} a_F\right)!\prod_{i\in T_v} k_i!},
%\]
%and optimizing $\prod_{v\in V(\Gamma)} d_v^{n(v)}$ (where $d_v := \sum_{F=(v,e)\in F_v} d_e$ satisfies $\sum_v d_v = 2D$) by Lagrange multipliers 
%shows that $\prod_{v} d_v^{n(v)}$ attains its maximum at $d_v = \frac{2D}{n+2\sharp E(\Gamma)} n(v)$, yielding 
%\begin{align*}
%\prod_{v\in V(\Gamma)} d_v^{n(v)} &\le \left(\frac{2D}{n+2\sharp E(\Gamma)}\right)^{n+2\sharp E(\Gamma)} \prod_{v\in V(\Gamma)} n(v)^{n(v)} \\
%&\le e^{2D} \prod_{v\in V(\Gamma)} \left(e^{n(v)} n(v)!\right) \le e^{n+4D} \prod_{v\in V(\Gamma)} n(v)!.
%\end{align*}
%\end{remark}

\subsection{Sums over graphs and proof of Theorem \ref{thm:main}}
Combining \eqref{eq:aut-order}, \eqref{eq:loc-sum}, and \eqref{eq:single-graph-bound}, we have uniformly for $a\in B_0$:
\begin{equation}\label{eq:correlator-reduced-to-graphs}
\left|\corr{\tau_{k_1}\gamma_{h_1},\dots,\tau_{k_n}\gamma_{h_n}}^T_{g,n,d}(a)
\right|
\le A_2(g,X)^{n+|d|+1+\sum_{i=1}^n k_i} 
\sum_{\hat{\Gamma}} \frac{1}{\sharp \Aut(\hat{\Gamma})} \prod_{v\in V(\Gamma)} n(v)!.
\end{equation}
Thus the proof of \eqref{eq:main-estimate} reduces to bounding the sum over decorated graphs $\hat{\Gamma}$ of genus $g$, degree $d$, and $n$ tails.

First, we reduce the sum over connected graphs $\Gamma$ with $h := \dim H^1(\Gamma) \in \{0,1,\dots,g\}$ to a sum over trees. 
For each decorated graph $\hat{\Gamma}$ with $\dim H^1(\Gamma)=h$, choose $h$ edges of $\Gamma$,  together with an orientation for each edge, whose removal breaks all cycles 
and turns $\Gamma$ into a tree $\tau$. 
Cutting these $h$ oriented edges creates $2h$ new tails, which we label by $n+1, n+2,\dots, n+2h$ so that the tail $n+2j-1$ is paired with $n+2j$ for each $1\le j\le h$. 
The decorations $p,l,d,g,m$ of $\hat{\Gamma}$ naturally induce decorations on $\tau$. Moreover, $\tau$ is equipped with additional decorations: for each pair $(n+2j-1,n+2j)$ of tails ($1\le j\le h$), we record the degree $d_j$ and the $T$-invariant rational curve $l_j$ originally associated with the corresponding cut edge. We let $\hat{\tau}$ denote the tree $\tau$ equipped with the decorations $(p,l,d,g,m, (d_j,l_j)_{1\le j\le h})$. 
Note that $\hat{\Gamma}$ is uniquely reconstructed from the resulting decorated tree $\hat{\tau}$ by gluing the pairs of tails 
$(n+2j-1, n+2j)$ for $1\le j\le h$.  Every automorphism of $\hat{\tau}$ induces an automorphism of $\hat{\Gamma}$ and hence $\sharp \Aut(\hat{\tau}) \le \sharp \Aut(\hat{\Gamma})$. 
Moreover, the number $n(v)$ at every vertex $v$ in $\tau$ is the same as that in $\Gamma$. 
Therefore, 
\begin{equation}\label{eq:cut-to-trees}
\sum_{\substack{\text{$\hat{\Gamma}$: genus $g$} \\ \text{degree $d$, $n$ tails}}} \frac{1}{\sharp \Aut(\hat{\Gamma})} \prod_{v\in V(\Gamma)} n(v)!
\le \sum_{h=0}^g \sum_{\substack{\text{$\hat{\tau}$: decorated tree} \\ \text{degree $d$, genus $g-h$} \\ \text{$n+2h$ tails}}} \frac{1}{\sharp \Aut(\hat{\tau})} \prod_{v\in V(\tau)} n(v)!,   
\end{equation}
where the degree of $\hat{\tau}$ is defined to be $\sum_{e\in E(\tau)} d_e [l_e] + \sum_{j=1}^h d_j [l_j]$ and the genus of $\hat{\tau}$ is defined to be $\sum_{v\in V(\tau)} g(v)$. 
%where $\tilde{\Gamma}$ runs over isomorphism classes of decorated trees with $n+2h$ labelled tails, 
%$j = \sharp E(\tilde{\Gamma}) = \sharp E(\Gamma) - h \le D \le C|d|$ tree edges (together with the data $(l_e, d_e)$ on the $h$ cut edges), 
%$\sum_{e\in E(\Gamma)} d_e = D$, and $\sum_{v} g(v) = g-h$.

Second, let $0\le h\le g$ and let us fix the isomorphism class of an unlabelled tree $\tau$ with $k$ edges and $m = n+2h$ unlabelled tails that underlies a decorated tree $\hat{\tau}$ of degree $d$ and genus $g-h$. 
Let $S(\tau)$ be the set of all decorations $(p,l,d,g,m,(d_j,l_j)_{1\le j\le h})$ on $\tau$ of degree $d$ and genus $g-h$. 
We have the estimate $k\le C|d|$ as in  \eqref{eq:edge-degree-bound}.   
The automorphism group $\Aut(\tau)$ of the unlabelled tree $\tau$ acts naturally on $S(\tau)$, and the orbits $\mO \subset S(\tau)$ 
are in bijection with isomorphism classes of decorated trees $\hat{\tau}$ over $\tau$, with orbit size $\sharp \mO = \sharp \Aut(\tau) / \sharp \Aut(\hat{\tau})$. 
Hence 
\[
\sum_{\hat{\tau} \text{ over } \tau} \frac{1}{\sharp \Aut(\hat{\tau})} = \frac{\sharp S(\tau)}{\sharp \Aut(\tau)}.
\]
Using constants $C_9,C_{10},C_{11}>1$ depending only on $X$, we bound $\sharp S(\tau)$ as follows: 
\begin{itemize}
\item The number of choices for $p\colon V(\tau)\to X^T$, $l\colon E(\tau) \to E(X)$ and $(l_j)_{j=1}^h$ is at most $(\sharp X^T)^{k+1}(\sharp E(X))^{k+h} \le C_9^{|d|+g+1}$.
\item The number of degree assignments $d\colon E(\tau) \to \Z_{>0}$, $(d_j)_{j=1}^h$ with $\sum d_e + \sum d_j = D \le C|d|$ is at most 
$\sum_{D=k+h}^{\lfloor C|d|\rfloor} \binom{D-1}{k+h-1} = \binom{\lfloor C|d|\rfloor}{k+h} \le 2^{C|d|} \le C_{10}^{|d|}$. 
%$\binom{D-1}{k+h-1} \le 2^D \le C_{10}^{|d|}$. %(and $1$ if $D=k+h=0$).
\item The number of genus assignments $g\colon V(\tau)\to \Z_{\ge 0}$ with $\sum_{v\in V(\tau)} g(v) = g-h$ is $\binom{k+g-h}{k} \le 2^{k+g} \le C_{11}^{|d|+g}$.
\item The number of tail labellings $m\colon \{1,\dots,n+2h\}\xrightarrow{\sim} T(\tau)$ is $(n+2h)!$.
\end{itemize}
Consequently, for $0\le h\le g$,
\begin{equation}\label{eq:tree-sum-separated}
\sum_{\substack{\text{$\hat{\tau}$: decorated tree} \\ \text{degree $d$, genus $g-h$} \\ \text{$n+2h$ tails}}} \frac{1}{\sharp \Aut(\hat{\tau})} \prod_{v\in V(\tau)} n(v)! 
\le C_{12}^{|d|+g+1} (n+2h)! \sum_{k=0}^{\lfloor C|d|\rfloor} c_{k, n+2h},
\end{equation}
for some constant $C_{12}>1$ depending only on $X$, where we define 
\[
c_{k,m} := \sum_{\substack{\text{$\tau$: unlabelled tree} \\ \text{$k$ edges, $m$ tails}}} \frac{1}{\sharp \Aut(\tau)} \prod_{v\in V(\tau)} n(v)!.
\]

\begin{lemma}\label{lem:tree-generating-function}
There exists an absolute constant $C_{13} > 1$ such that for all $k\ge 0$ and $m\ge 0$,
\[
c_{k,m} \le C_{13}^{k+m+1}.
\]
\end{lemma}
\begin{proof}
For any tree $\tau$ with $k$ edges and $m$ tails and a vertex $v_0\in V(\tau)$, since $\prod_{v\in V(\tau)\setminus \{v_0\}} n(v) \le e^{\sum_{v\in V(\tau)\setminus \{v_0\}} n(v)} \le e^{2k+m}$, we have 
\[
\prod_{v\in V(\tau)} n(v)! \le e^{2k+m} n(v_0)! \prod_{v\in V(\tau)\setminus \{v_0\}} (n(v)-1)!. 
\]
Therefore $c_{k,m} \le e^{2k+m} \tilde{c}_{k,m}$ with 
\[
\tilde{c}_{k,m} := \sum_{(\tau, v_0)} \frac{1}{\sharp \Aut(\tau, v_0)} n(v_0)! \prod_{v\in V(\tau)\setminus \{v_0\}} (n(v)-1)!
\]
being the sum over isomorphism classes of \emph{rooted} trees $(\tau, v_0)$ with $k$ edges, $m$ unlabelled tails, and root vertex $v_0$. 

Consider the generating function $F(q,t) := \sum_{k,m\ge 0} \tilde{c}_{k,m} q^k t^m$. 
Removing the root vertex $v_0$ decomposes a rooted tree $(\tau, v_0)$ into $i_0$ tails attached to $v_0$ 
and $j$ rooted subtrees $(\tau_1, v_1),\dots,(\tau_j, v_j)$ joined to $v_0$ by edges,  
where each $v_i$ is adjacent to $v_0$ in $\tau$, $n_{\tau}(v_0) = i_0 + j$ and $n_{\tau}(v_i)  = n_{\tau_i}(v_i)+1$.  
Considering the $j!$ symmetries permuting the branches $(\tau_1, v_1),\dots,(\tau_j,v_j)$ and the $i_0!$ symmetries permuting the tails attached to $v_0$, we obtain the recursion
\[
\tilde{c}_{k,m} = \sum_{j\ge 0} \sum_{\substack{i_0 + m_1 + \cdots + m_j = m \\ j + k_1 + \cdots + k_j= k}}
\frac{(i_0+j)!}{i_0!\, j!} \tilde{c}_{k_1, m_1} \cdots \tilde{c}_{k_j, m_j}.
\]
In terms of the generating function $F(q,t)$, this recursion reads
\[
F(q,t) = \sum_{j\ge 0} \sum_{i_0\ge 0} \frac{(i_0+j)!}{i_0!\, j!} t^{i_0} (q F(q,t))^j 
= \sum_{N=0}^\infty (t + q F(q,t))^N = \frac{1}{1 - (t + q F(q,t))}.
\]
Thus $F(q,t)$ satisfies the quadratic equation $q F^2 - (1-t)F + 1 = 0$ with $F(0,0)=1$, so 
\[
F(q,t) = \frac{1 - t - \sqrt{(1-t)^2 - 4q}}{2q} = \frac{2}{1 - t + \sqrt{(1-t)^2 - 4q}}.
\]
Since $F(q,t)$ is holomorphic in a neighborhood of $(q,t)=(0,0)$, its Taylor coefficients $\tilde{c}_{k,m}$ grow at most exponentially in $k+m$, proving the lemma.
\end{proof}

Using Lemma \ref{lem:tree-generating-function}, \eqref{eq:tree-sum-separated}, \eqref{eq:cut-to-trees} and $(n+2h)! \le n! (2h)!\, 2^{n+2h}\le n! (2g)! 2^{n+2g}$, we obtain 
\[
\sum_{\hat{\Gamma}} \frac{1}{\sharp \Aut(\hat{\Gamma})} \prod_{v\in V(\Gamma)} n(v)!
\le C_{14}^{n+|d|+g+1} n! (2g)!,
\]
for some $C_{14}>1$. This together with \eqref{eq:correlator-reduced-to-graphs} yields the uniform bound on $B_0$:
\begin{equation}\label{eq:uniform-ball-bound}
\sup_{a\in B_0} \left|\corr{\tau_{k_1}\gamma_{h_1},\dots,\tau_{k_n}\gamma_{h_n}}^T_{g,n,d}(a)\right|
\le A_3(g,X)^{n+|d|+1+\sum_{i=1}^n k_i} n! 
\end{equation}
for some $A_3(g,X)>1$ depending only on $g$ and $X$.

We now complete the proof of Theorem \ref{thm:main}. 
For each tuple $(h_1,\dots,h_n)$ and $(k_1,\dots,k_n)$, the equivariant correlator 
\[
P(a) := \corr{\tau_{k_1}\gamma_{h_1},\dots,\tau_{k_n}\gamma_{h_n}}^T_{g,n,d}(a) 
\]
is a homogeneous polynomial in $a\in\mathfrak{t}_\C$ of degree 
\begin{equation}\label{eq:equiv-degree-m}
m = \sum_{i=1}^n \left(k_i + \frac{1}{2}\deg \gamma_{h_i}\right) - ((1-g)(s-3) + n + \langle c_1(X), d\rangle). 
\end{equation}
Note that $P=0$ if $m<0$. 
Note that \eqref{eq:equiv-degree-m} implies 
\begin{equation}
\label{eq:sum-ki-bound}
\sum_{i=1}^n k_i \le m + |(1-g)(s-3)| + n + |\langle c_1(X), d\rangle| \le m + C_{15}(n+|d|+g+1)
\end{equation}
for some $C_{15}>1$ depending only on $X$. 
Recall that $B_0\subset \mathfrak{t}_\C$ is a ball of radius $\delta>0$ (Lemma \ref{lem:choice-of-ball}); let $a_0\in B_0$ denote the centre. For $u\in \mathfrak{t}_\C$, consider the univariate polynomial $f(z) = P(a_0+z u)$ in $z$. The homogeneity of $P$ implies 
\[
f(z)  = z^m P(u) + O(z^{m-1}), 
\]
and hence the Cauchy integral formula gives the estimate
\[
|P(u)| = \left| \frac{1}{2\pi \sqrt{-1}} \int_{|z|=1} f(z)\frac{dz}{z^{m+1}} \right| 
\le \sup_{a\in B_0} |P(a)| 
\]
for any $u\in \mathfrak{t}_\C$ with $|u|\le \delta$. 
Thus by \eqref{eq:uniform-ball-bound} and the homogeneity of $P$ again, for any $r>0$ and any $u\in \mathfrak{t}_\C$ with $|u|\le r$, we have
\begin{align}
\label{eq:estimate-P} 
\begin{split}  
|P(u)| 
& = (r/\delta)^m |P(\tfrac{\delta}{r} u)| 
\le (r/\delta)^m A_3(g,X)^{n+|d|+1+\sum k_i}  n! \\
& \le (r A_3(g,X)/\delta)^m A_4(g,X)^{n+|d|+1} n! 
\end{split} 
\end{align} 
for some $A_4(g,X)>1$ depending only on $g$ and $X$, 
where we used \eqref{eq:sum-ki-bound} in the second line. Theorem \ref{thm:main}(1) follows immediately by setting $r=\delta/A_3(g,X)$. 

When $|a|\le r$, $|t_{k,i}|\le r$ and $|Q^d|\le r^{\langle \omega, d\rangle}$ for all effective classes $d\in H_2(X,\Z)$, using \eqref{eq:estimate-P}, we have 
\begin{align*} 
|\mF_g(\bt,a;Q)|
& \le \sum_{n\ge 0} \sum_{d\in \Eff(X)} \frac{1}{n!}\sum_{\substack{k_1,\dots,k_n\ge 0 \\ 1\le h_1,\dots,h_n\le N}}
|\corr{\tau_{k_1}\gamma_{h_1},\dots,\tau_{k_n}\gamma_{h_n}}^T_{g,n,d}| 
|t_{k_1,h_1}|\cdots |t_{k_n,h_n}| 
|Q^d|\\ 
& \le \sum_{n\ge 0} \sum_{d\in \Eff(X)}\sum_{\substack{k_1,\dots,k_n\ge 0 \\ 1\le h_1,\dots,h_n\le N \\
\text{such that $m\ge 0$}} }
 (r A_3(g,X)/\delta)^m A_4(g,X)^{n+|d|+1} r^{n+\langle \omega,d\rangle} 
\end{align*} 
where $m$ in the second line is given by \eqref{eq:equiv-degree-m} and $\Eff(X) \subset H_2(X,\Z)$ denotes the set of effective classes. 
Assuming now $0<r\le \delta/(2 A_3(g,X))<1$, we have
\begin{align*} 
|\mF_g(\bt,a;Q)| & \le \sum_{n\ge 0} \sum_{d\in \Eff(X)}
\sum_{\substack{k_1,\dots,k_n\ge 0 \\ 1\le h_1,\dots,h_n\le N \\ \text{such that $m\ge 0$}}}
 2^{-m} A_4(g,X)^{n+|d|+1} r^{n+\langle \omega,d\rangle} \\
 & \le \sum_{n\ge 0} \sum_{d\in \Eff(X)} N^n \sum_{k_1,\dots,k_n\ge 0} 
 2^{-(\sum k_i) + (n+|d|+g+1)C_{15}} A_4(g,X)^{n+|d|+1} r^{n+\langle \omega,d\rangle} \\ 
 & = \sum_{n\ge 0} \sum_{d\in \Eff(X)} 
2^{(g+1)C_{15}}A_4(g,X) \left(2N 2^{C_{15}}A_4(g,X)r\right)^n \left(2^{C_{15}} A_4(g,X) r^{\varepsilon}\right)^{|d|} 
\end{align*}  
where we used \eqref{eq:sum-ki-bound} in the second line and $\varepsilon>0$ is a constant satisfying $\langle \omega,d\rangle\ge \varepsilon |d|$ for all $d\in \Eff(X)$. This converges if both $2N2^{C_{15}}A_4(g,X) r$ and $2^{C_{15}} A_4(g,X) r^\varepsilon$ are less than $1$ (note that the number of classes $d\in H_2(X,\Z)$ with $|d|\le R$ grows polynomially in $R$). We have established Theorem \ref{thm:main}(2).

\end{document}